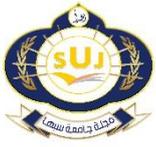
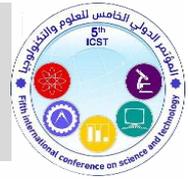

# Roughness in Anti Semigroup

Faraj. A. Abdunabi, Ahmed shletiet, *Najah. A. Bosaif

Department of Mathematics, Faculty of Science, AJdabyia University, Libya

**A B S T R A C T**
In this paper, we present the concepts of the upper and lower approximations of Anti-rough subgroups, Anti-rough subsemigroups, and homeomorphisms of Anti-Rough anti-semigroups in approximation spaces. Specify the concepts of rough in Finite anti-groups of types (4) are studies. Moreover, some properties of approximations and these algebraic structures are introduced. In addition, we give the definition of homomorphism anti-group.

## الغموض في شبه الزمرة المضادة

فرج أرخيص عبدالنبي و أحمد أبراهيم المبروك أشليتيتب و *نجاح عبدالقادر بوسيف

قسم الرياضيات، كلية العلوم، جامعة جدابيا، أجدابيا، ليبيا

**الملخص**
في هذا البحث، نقدم مفاهيم التقريب العلوي والسفلي لزمر جزئية المضادة، والزمر شبه الجزئية المضادة، والتشاكل للزمر المضادة والشبه مضادة في فضاء التقريبات. كما سندرس بشكل خاص مفاهيم الضبابية في شبه الزمر المضادة المحدودة من الأنواع (4). علاوة على ذلك، نقدم بعض خصائص التقريبات وهذه الهياكل الجبرية. ونعطي تعريف التشاكل anti-group.

## 1-Introduction

Pawlak [1] 1982 introduced the concept of the rough set theory as a new and good tool for modeling in an information system. This theory has prompted many types of interest by many researchers. It has developed amazingly in pure mathematics. Some authors have studied the algebraic structures of rough sets such as Bonikowaski [2], Iwinski [3], and Pomykala and Pomykala [4]. Miao et al.[5] have improved the rough group and rough subgroup and considered some properties. In 1994, Biswas and Nanda [6] introduced the definition of a rough group depending on the upper approximate, not on the lower approximation. B.Davvaz in [7], studied the concept of rough subring with respect to an ideal. Yao in [8] considered the concepts of lower and upper approximations on the lattice. In addition, some properties of the lower and the upper approximations with respect to the normal subgroups studied in [9]. The concepts of rough set theory build on lower and upper approximations. The upper approximation of a given set is the union of all the equivalence classes that are subsets of the set, and the upper approximation is the union of all the equivalence classes that are intersection with a non-empty set. The main purpose of this paper is to introduce rough anti-semigroups of Finite anti-groups of types (4). In addition, some properties of approximations of these algebraic structures are introduced. Moreover, the notion of Anti-Rough semigroups was introduced.

However, our definition of rough anti-semigroup is similar to the definition of rough groups.

## 2-Preliminaries

In this section, the most important concepts of rough set theory needed for this research are presented.

Suppose that $\sim$ an equivalence relation on a universe set $U$ ($\emptyset, finite$). The pair $(U, \sim)$ is called an approximation space. The family of all equivalent classes $[x_\sim]$ denotes by $U/\sim$. For any $\backslash M \subseteq U$, write $M^c$ to denote the complementation of $M$ in U.

**Definition 2.1**: Let $(U, \sim)$ be an approximation space. Define the upper approximation of $M$ by $\overline{M} = \{x \in U: [x]_\sim \cap M \neq \emptyset\}$ and the lower approximation of $M$ by

$\sim M = \{x \in U: [x]_\sim \subseteq M\}$. The difference $BM_\sim = \overline{M} - \sim M$ is called the boundary. If $BM_\sim = \emptyset$, we say $M$ is an exact (crisp) set otherwise, $M$ is a Rough set (inexact).

**Preposition 2-1**: Let $(U, \sim)$ be an approximation space and $X, Y \subseteq U$, we have:

1) $\sim X \subseteq X \subseteq \overline{X}$;
2) $\sim \emptyset = \overline{\sim \emptyset}, \sim U = \overline{U}$,
3) $\sim (X \cup Y) \supseteq \sim (X) \cup \sim (Y)$,

Corresponding author:
E-mail addresses: Najahboseaf@hotmail.com, (F. A. Abdunabi) faraj.a.abdunabi@uoa.edu.ly, (A. shletiet) Ahmed.Shlitite@uoa.edu.ly

4) $\sim(X \cap Y) = \sim(X) \cap \sim(Y)$,
5) $\overline{\sim(X \cup Y)} = \overline{\sim(X)} \cup \overline{\sim(Y)}$.
6) $\overline{\sim(X \cap Y)} \subseteq \overline{\sim(X)} \cap \overline{\sim(Y)}$.
7) $\overline{X^c} = (\underline{\sim X})^c \cdot \underline{\sim X^c} = (\overline{\sim X})^c$.
8) $\sim(\sim X) = \overline{\sim(\sim X)} = \underline{\sim X}$. –
9) $\overline{(\sim(\overline{\sim X})} = \sim(\overline{\sim X}) = \overline{\sim X}$.

**Proposition 2-2** [8] Let $(U, R)$ be an approximation space. Let $X$ and $Y$ be nonempty subsets of $U$. Then
1) $\overline{\sim X} \ \overline{\sim Y} = \overline{\sim XY}$.
2) $\underline{\sim X} \ \underline{\sim Y} \subseteq \underline{\sim XY}$.

**Definition 2.2**[10]. Suppose that G is a nonempty set.
Let $* : \mathcal{R} \times \mathcal{R} \to \mathcal{R}$ be binary operations defined on $G$. The $(G,*.)$ is called a group if satisfy the following conditions:
C1: For all x, y ∈ G, x*y ∈ G ;
C2: For all x, y, z ∈ G , $x * (y * z) = (x * y) * z$ ;
C3: For all x ∈ G , there exists e ∈ G such that $x * e = e * x = x$;
C4: For all x ∈ G, there exists −x ∈ G such that $x * (-x) = (-x) * x = e$;
If we have,
C5: For all x, y ∈ G, $x * y = y * x$, then (G, *) is called a commutative group.

**Definition 2.3**[10]. A semigroup S is an algebraic structure on a nonempty set together with an associative binary operation. That means, a semigroup is a set together with a binary operation "*" that satisfies C1,C2.

**Definition 2.4**. A nonempty subset $H$ of a semigroup S is said to be a subsemigroup of S, if $a * b \in H \ for \ all \ a, b \in S$.

**Definition 2.5**. An anti-group $\mathfrak{C}$ is an alternative to the group G that has at least one anti-Law or at least one flowing conditions:
For all the duplets $(x, y) \in \mathfrak{C}, x * y \notin \mathfrak{C}$ ;
C7: For all the triplets $(x, y, z) \in \mathfrak{C}, x * (y * z) \neq (x * y) * z$;.
C8: There does not exist an element e ∈ $\mathfrak{C}$ such that $x * e = e * x = x \ \forall x \in \mathfrak{C}$.
C9: There does not exist u ∈$\mathfrak{C}$ such that $* u = u * x = e \ \forall x \in \mathfrak{C}$.

**Definition 2.6**. An anti- abelian-group $\mathfrak{C}$ is an alternative to the classical an abelian group G that has at least one Anti-Law or at least one of {C6, C7,C8, C9} and
C10: For all the duplets $(x, y) \in \mathfrak{C}, x * y \neq y * x$.
A particular class of Anti-groups ($\mathfrak{C}$, *) where G4 is totally false for all the elements of $\mathfrak{C}$ while C1, C2, C3 and C5 are either partially true, partially indeterminate or partially false for some elements of $\mathfrak{C}$.

**Proposition 2.3**. Let ($\mathfrak{C}$, *) be an Anti- group of type-AG(4) and let $g, x, y \in \mathfrak{C}$. Then
1) $g * x = g * y \nRightarrow x = y$.
2) $x * g = y * g \nRightarrow x = y$.

**Definition 2.9**. Let ($\mathfrak{C}$, *) be an anti-group of type-AG(4) and let $A$ and $B$ be an anti-Subgroups of $\mathfrak{C}$. The set A*B is defined by $A * B = \{x \in : x = h * k \ for \ some \ h \in A, k \in B\}$.

**3-Roughness in Anti- semigroups**
In this section, the notions of rough anti-semigroup and rough sub semigroup on an approximation space are introduce and study some of its properties.

**Definition3-1**.[8] Suppose that $(U, \sim)$ is an approximation space and (*) be a binary operation defined on U. A subset $A$ of U is called a rough anti- semigroup on approximation space, provided the following properties are satisfied:
1) For all $x, y \in A, x * y \in \overline{\sim A}$,

2) For all $x, y, z \in A, (x * y) * z = x * (y * z)$ property holds in $\underline{\sim A}$.

**Example 3.1.** Let U = {1, 2, 3, 4, 5, 6} be a universe of discourse and $\mathfrak{C}$ = {1, 2,3,5} be a subset of U.
Let * be a binary operation defined on $\mathfrak{C}$ as shown in the Cayley table below

| * | 1 | 2 | 3 | 5 |
|---|---|---|---|---|
| 1 | 4 | 1 | 3 | 5 |
| 2 | 1 | 4 | 5 | 3 |
| 3 | 2 | 1 | 6 | 5 |
| 5 | 1 | 2 | 3 | 6 |

It is evident from the above table that C1, C2, C3, C5 are either partially true or partially false with respect to * but C4 is false for all the elements of $\mathfrak{C}$. Hence ($\mathfrak{C}$, *) is a finite Anti-group of $\mathfrak{C}$. A classification of U is $U/\sim = \{E1, E2, E3\}, where \ E1 = \{1, 2, 3\}, E2 = \{4, \}, E3 = \{5\}$.
let $A$ = {1, 2, 5}, Let * be defined on $A$ as shown in the Cayley tables below:

| * | 1 | 2 | 5 |
|---|---|---|---|
| 1 | 4 | 1 | 5 |
| 2 | 1 | 4 | 3 |
| 5 | 1 | 2 | 6 |

(e-6)
It can easily be seen from the tables that $A$ is an anti-Subgroup of type-AG(4). $\overline{\sim A}$ = {1, 2, 3, 4} . From Definition 3-1, $A \subseteq$ U is a rough anti-semigroup.

**Definition 3-2.** Suppose that (U, ~) be an approximation space and (*) be a binary operation defined on U. Let $A$ be a rough Anti-semigroup and H a nonempty subset of $A$. A nonempty subset H of a rough anti-semigroup $A$ is said to be a rough anti-subsemigroup of $A$, if $a * b \in \overline{\sim H}$ for all a, b ∈ H, i.e., $HH \subseteq \overline{\sim H}$.

**Example 3.2.** Consider example 3-1. Let $B$= ($\mathfrak{C}$, *) be the $B$= {2,3, 5} a subset of $\mathfrak{C}$ and * be defined on $B$ as shown in the Cayley tables below:

| * | 2 | 3 | 5 |
|---|---|---|---|
| 2 | 4 | 5 | 3 |
| 3 | 1 | 6 | 5 |
| 5 | 2 | 3 | 6 |

.
It can easily be seen from the tables that is an anti-Subgroup of $\mathfrak{C}$. $\overline{\sim B}$ = {1, 2, 3, 5} . From Definition 3-1, $B \subseteq$ U. is a rough anti-semigroup.

**Proposition 3-1**. Suppose that $(U, \sim)$ be an approximation space and (*) be a binary operation defined on U. Suppose that $A$ and $B$ be two rough anti sub semigroups of the rough anti-semigroup $A$. Then $\overline{\sim(A)} \cap \overline{\sim(B)} \subseteq \overline{\sim(A \cap B)}$.
A sufficient condition for intersection of two rough anti-sub semigroups of a rough anti-semigroup be a rough anti subsemigroup is $\overline{\sim(A)} \cap \overline{\sim(B)} = \overline{\sim(A \cap B)}$.

**Example 3.3.** Consider example 3.1 and 3. 2. $A$ = {1, 2, 5} and = {2,3, 5}, then $A \cap B$ ={2,5}
then $\overline{\sim A}$ = {1, 2, 3, 4}∩ $\overline{\sim B}$ = {1, 2, 3, 5}= {{1, 2, 3}, $\overline{\sim(A \cap B)}$ = {1, 2, 3, 5}.

**4-Homomorphism of rough Anti-Group**
Suppose that ($\mathfrak{C}$, *) and ($\mathfrak{B}$, ∘) be any two anti-groups of type-AG(4). The mapping $\varphi : \mathfrak{C} \to \mathfrak{B}$ is called an Anti-group Homomorphism if φ does not preserve the binary operations * and ∘ that is for all the duplet $(x, y) \in \mathfrak{C}$, we have $\varphi(x * y) \neq \varphi(x) \circ \varphi(y)$.

The kernel of φ denoted by Kerφ is defined by

Kerφ = {x : φ(x) = e$\mathfrak{B}$ for at least one e$\mathfrak{B} \in \mathfrak{B}$} where e$\mathfrak{B}$ is a NeutroNeutral Element in $\mathfrak{B}$. The image of φ denoted by Imφ is defined by Imφ = {y ∈ : y = φ(x) for some x ∈ $\mathfrak{C}$}.

If in addition φ is an anti bijection, then φ is called an Anti-group Isomorphism.

Suppose that Let $(U1, \sim), (U2, \rho)$ be two approximation spaces, and $(\cdot)$ be binary operation over universes U1 and ,$(\circ)$ over universes U2

**Definition 4.1.** Let *A*⊂U1 and *B*⊂U2 be rough anti-semigroups. If there exists a surjection $\phi : \overline{\sim(A)} \to \overline{\sim(B)}$ such that $\phi(x \cdot y) = \phi(x) \circ \phi(y)$ for all x, y ∈ $\overline{\sim A}$ then ϕ is called a rough homomorphism and *A*, *B* are called rough homomorphic semigroups.

**Definition 4.2** Let $\mathfrak{C} \subset$ U1, $\mathfrak{B} \subset$ U2 be rough anti groups. If there exists a surjection ϕ : $\overline{\sim(\mathfrak{C})} \to \overline{\sim(\mathfrak{B})}$ such that ϕ(x · y) = ϕ(y) ∘ ϕ(x) for all x, y ∈ $\overline{\sim\mathfrak{C}}$ then ϕ is called a rough anti homomorphism.

**Proposition 4.1.** Let $\mathfrak{C}$ be a rough anti-group and φ1 be a rough anti-homomorphism and φ2 be a rough homomorphism on $\mathfrak{C}$. Then the composition φ1οφ2 is a rough anti-homomorphism on $\mathfrak{C}$.

Proof. Let $\mathfrak{C}$ be a rough anti-group and let φ1 be a rough anti-homomorphism on $\mathfrak{C}$ and φ2 be a rough homomorphism on $\mathfrak{C}$. Then φ1, φ2 : : $\overline{\sim(\mathfrak{C})} \to \overline{\sim(\mathfrak{B})}$ such that ∀x, y ∈ : $\overline{\sim(\mathfrak{C})}$, φ1(x ∗ y) = φ1(y) ∗ φ1(x) and φ(x ∗ y) = φ2(x) ∗ φ2(y) Now ∀ x, y ∈: $\overline{\sim(\mathfrak{C})}$ (φ1οφ2)(x ∗ y) = φ1(φ2(x ∗ y)) = φ1(φ2(x) ∗ φ2(y)) = (φ1οφ2)(y) ∗ (φ1οφ2)(x) Therefore, φ1οφ2 is a rough anti-homomorphism on $\mathfrak{C}$.

**Proposition 4.2.** Let $\mathfrak{C}$ be a rough anti-group and φ1 and φ2 be two rough anti-homomorphisms on $\mathfrak{C}$. Then the composition φ1οφ2 is a rough homomorphism on $\mathfrak{C}$.

Proof. Let $\mathfrak{C}$ be a rough anti-group and let φ1, φ2 be two rough anti-homomorphisms on $\mathfrak{C}$. Then φ1, φ2 : $\overline{\sim(\mathfrak{C})} \to \overline{\sim(\mathfrak{C})}$ such that ∀x, y ∈ $\overline{\sim(\mathfrak{C})}$ φ1(x ∗ y) = φ1(y) ∗ φ1(x) and φ2(x ∗ y) = φ2(y) ∗ φ2(x). Now ∀x, y ∈$\overline{\sim(\mathfrak{C})}$ (φ1οφ2)(x ∗ y) = φ1(φ2(x ∗ y)) = φ1(φ2(y) ∗ φ2(x)) = (φ1οφ2)(x) ∗ (φ1οφ2)(y) Therefore, φ1οφ2 is a rough homomorphism on $\mathfrak{C}$ .

**Conclusion**

The concepts of rough in Finite anti-groups of types (4) introduced in this paper. Moreover, some properties of approximations of these algebraic structures are studies and considers. However, the definition of homomorphism anti-group is given.

Acknowledgments: The private discussions and suggestions of staff of mathematics department of Ajdabyia University for help and suggests on this paper.
The valuable comments and suggestions of all the anonymous reviewers are equally acknowledged.